\documentclass[12pt]{amsart}
\numberwithin{equation}{section}

\usepackage{amssymb}
\def\cb{{\mathcal B}}

\def\ch{{\mathcal H}}

\def\ck{{\mathcal K}}

\def\bt{{\mathbb T}}

\def\a{\alpha}
\def\b{\beta}
 \def\G{\Gamma}
\def\d{\delta}  \def\D{\Delta}

\def\eps{\varepsilon}

\def\m{\mu}

\def\s{\sigma}

\def\z{\zeta}

\newtheorem{thm}{Theorem}[section]
\newtheorem{lem}[thm]{Lemma}

\newtheorem{prop}[thm]{Proposition}

\def\di{\mathop{\rm d}\!}

\newcommand{\nn}{\nonumber}

\begin{document}

\title[entangled ergodic theorem]
{on the entangled ergodic theorem}
\author{Francesco Fidaleo}
\address{Francesco Fidaleo,
Dipartimento di Matematica,
Universit\`{a} di Roma Tor Vergata, 
Via della Ricerca Scientifica 1, Roma 00133, Italy} \email{{\tt
fidaleo@mat.uniroma2.it}}

\date{December 11, 2005}

\begin{abstract}
Let $U$ be a unitary operator acting on the Hilbert space $\ch$, and 
$\a:\{1,\dots, m\}\mapsto\{1,\dots, k\}$ a partition of the 
set $\{1,\dots, m\}$.
We show that the ergodic average 
$$
\frac{1}{N^{k}}\sum_{n_{1},\dots,n_{k}=0}^{N-1}
U^{n_{\a(1)}}A_{1}U^{n_{\a(2)}}\cdots 
U^{n_{\a(m-1)}}A_{m-1}U^{n_{\a(m)}}
$$
converges in the weak 
operator topology if the $A_{j}$ belong to the algebra of 
all the compact operators on $\ch$. We write esplicitely the 
formula for these ergodic averages in the case of pair--partitions. 
Some results without any restriction on the operators $A_{j}$
are also presented in the almost periodic case.

\end{abstract}

\maketitle

\section{introduction}

An entangled ergodic theorem was formulated in \cite{AHO} in 
connection with the quantum central limit theorem. Namely, let $U$ be 
a unitary operator on the Hilbert space $\ch$, and for $m\geq k$, 
$\a:\{1,\dots,m\}\mapsto\{1,\dots, k\}$ a partition of the 
set $\{1,\dots,m\}$.\footnote{A partition $\a:\{1,\dots,m\}\mapsto\{1,\dots, k\}$
of the set made of $m$ elements in $k$ parts is nothing but a surjective 
map, the parts of $\{1,\dots,m\}$ being the preimages 
$\{\a^{-1}(\{j\})\}_{j=1}^{k}$.}
The {\it entangled ergodic theorem} concerns 
the convergence in the strong, or merely weak (s--limit, or w--limit 
for short) operator topology, of the multiple Cesaro mean
\begin{equation}
\label{0}
\frac{1}{N^{k}}\sum_{n_{1},\dots,n_{k}=0}^{N-1}
U^{n_{\a(1)}}A_{1}U^{n_{\a(2)}}\cdots 
U^{n_{\a(2k-1)}}A_{m-1}U^{n_{\a(m)}}
\end{equation}
$A_{1},\dots,A_{m-1}$ being bounded operators acting on $\ch$.

Expressions like \eqref{0} naturally appear also in \cite{NSZ}, in the study of the 
multiple mixing. The entangled ergodic 
theorem is a generalization of the well--known mean ergodic theorem 
due to von Neumann (cf. \cite{RS}) 
$$
\mathop{\rm s\!-\!lim}_{N}\frac{1}{N}\sum_{n=0}^{N-1}U^{n}=E_{1}\,,
$$
$E_{1}$ being the selfadjoint projection onto the eigenspace of the 
invariant vectors for $U$.

Therefore, it is natural to address 
the systematic investigation of the conditions under which the entangled ergodic theorem
holds true. A first attempt was done in \cite{L}, where some facts 
concerning the structure 
of the above ergodic average were pointed out. Unfortunately, there 
is yet no general result on the entangled ergodic theorem.

In the present note, we show that the ergodic average 
\eqref{0} converges in the weak 
operator topology if the $A_{k}$ belong to $\ck(\ch)$, the algebra of 
all the compact operators acting on $\ch$. We write down the 
formula for those ergodic averages in the case of pair--partitions as
\begin{align*}
&\mathop{\rm w\!-\!lim}_{N}\bigg\{\frac{1}{N^{k}}\sum_{n_{1},\dots,n_{k}=0}^{N-1}
U^{n_{\a(1)}}A_{1}U^{n_{\a(2)}}\cdots 
U^{n_{\a(2k-1)}}A_{2k-1}U^{n_{\a(2k)}}\bigg\}\\
=&\sum_{z_{1},\dots,z_{k}\in\s_{\mathop{\rm pp}}^{\mathop{\rm a}}(U)}
E_{z^{\#}_{\a(1)}}A_{1}E_{z^{\#}_{\a(2)}}\cdots
E_{z^{\#}_{\a(2k-1)}}A_{2k-1}E_{z^{\#}_{\a(2k)}}
\end{align*}
(see below for the precise definition of the r.h.s. of this formula), and conjecture 
that it holds true for each set $A_{1},\dots,A_{2k-1}$ of bounded operators. 

Finally, we present some results on the entangled ergodic theorem relative to the case 
when the 
dynamics induced by the unitary $U$ on $\ch$ is almost periodic, that 
is when $\ch$ is generated by the eigenvectors of $U$ (cf. \cite{NSZ}) 
without any restriction on the operators $A_{i}$.

\section{the entangled ergodic theorem for compact operators}

Let $U\in\cb(\ch)$ be a unitary operator, and 
$\a:\{1,\dots, 2k\}\mapsto\{1,\dots, k\}$ a pair--partition of the 
set $\{1,\dots, 2k\}$. Define
\begin{equation*}
\s_{\mathop{\rm pp}}^{\mathop{\rm a}}(U):=
\big\{z\in\s_{\mathop{\rm pp}}(U)\,\big|\,zw=1\,\text{for 
some}\,w\in\s_{\mathop{\rm pp}}(U)\big\}
\end{equation*}
where $\s_{\mathop{\rm pp}}(U)=
\big\{z\in\bt\,\big|\,z\,\text{is an eigenvalue of}\,U\big\}$, see 
e.g. \cite{RS}.

Consider, for each finite subset 
$F\in\s_{\mathop{\rm pp}}^{\mathop{\rm a}}(U)$ and 
$\{A_{1},\dots,A_{2k-1}\}\subset\cb(\ch)$, the following 
operator
\begin{equation}
\label{opses}
S^{F}_{\a;A_{1},\dots,A_{2k-1}}:=
\sum_{z_{1},\dots,z_{k}\in F}E_{z^{\#}_{\a(1)}}A_{1}E_{z^{\#}_{\a(2)}}\cdots
E_{z^{\#}_{\a(2k-1)}}A_{2k-1}E_{z^{\#}_{\a(2k)}}
\end{equation}
together with the sesquilinear form
\begin{equation*}
s^{F}_{\a;A_{1},\dots,A_{2k-1}}(x,y):=\big\langle 
S^{F}_{\a;A_{1},\dots,A_{2k-1}}x,y\big\rangle\,,
\end{equation*}
where the pairs $z^{\#}_{\a(i)}$ are alternatively $z_{j}$ and
$\bar z_{j}$ whenever $\a(i)=j$, and $E_{z}$ is the selfadjoint 
projection on the eigenspace corresponding to the eigenvalue 
$z\in\s_{\mathop{\rm pp}}(U)$.\footnote{If for example, $\a$ is the pair--partition
$\{1,2,1,2\}$ of four elements, we write
$$
S^{F}_{\a;A,B,C}=
\sum_{z,w\in F}E_{z}AE_{w}BE_{\bar z}CE_{\bar w}\,.
$$

If we have the trivial pair--partition of two elements, we write
${\displaystyle S^{F}_{A}=
\sum_{z\in F}E_{z}AE_{\bar z}}$,
and
${\displaystyle S_{A}=
\sum_{z\in\s_{\mathop{\rm pp}}^{\mathop{\rm a}}(U)}E_{z}AE_{\bar z}}$
for its limit in the strong operator topology (cf. Lemma 
\ref{2}), omitting the subscript $\a$.}
\begin{lem}
\label{1}
We have for the above sesquilinear form,
$$
\big|s^{F}_{\a;A_{1},\dots,A_{2k-1}}(x,y)\big|\leq\|x\|\|y\|\prod_{j=1}^{2k-1}\|A_{j}\|\,,
$$
uniformly for $F$ finite subsets of $\s_{\mathop{\rm pp}}^{\mathop{\rm a}}(U)$.

\end{lem}
\begin{proof}
To simplify matter, we deal with a particular case. The computation can be 
easily generalized to all the situations. Consider for example the 
pair--partition $\a$ of six elements given by $\{1,2,1,3,2,3\}$. We get
\begin{align*}
\big|s^{F}_{\a;A_{1},\dots,A_{5}}&(x,y)\big|^{2}
\leq\|y\|^{2}\big\|\sum_{z\in F}E_{z}A_{1}\sum_{w\in F}
E_{w}A_{2}E_{\bar z}A_{3}\sum_{\z\in F}E_{\z}A_{4}
E_{\bar w}A_{5}E_{\bar\z}x\big\|^{2}\\
=&\|y\|^{2}\sum_{z\in F}\big\|E_{z}A_{1}\sum_{w\in F}
E_{w}A_{2}E_{\bar z}A_{3}\sum_{\z\in F}E_{\z}A_{4}
E_{\bar w}A_{5}E_{\bar\z}x\big\|^{2}\\
\leq&(\|y\|\|A_{1}\|)^{2}\sum_{z\in F}\big\|\sum_{w\in F}
E_{w}A_{2}E_{\bar z}A_{3}\sum_{\z\in F}E_{\z}A_{4}
E_{\bar w}A_{5}E_{\bar\z}x\big\|^{2}\\
=&(\|y\|\|A_{1}\|)^{2}\sum_{z,w\in F}
\big\|E_{w}A_{2}E_{\bar z}A_{3}\sum_{\z\in F}E_{\z}A_{4}
E_{\bar w}A_{5}E_{\bar\z}x\big\|^{2}\\
\leq&(\|y\|\|A_{1}\|\|A_{2}\|)^{2}\sum_{w\in F}\sum_{z\in F}
\big\|E_{\bar z}A_{3}\sum_{\z\in F}E_{\z}A_{4}
E_{\bar w}A_{5}E_{\bar\z}x\big\|^{2}\\
\leq&(\|y\|\prod_{j=1}^{3}\|A_{j}\|)^{2}\sum_{w\in F}
\big\|\sum_{\z\in F}E_{\z}A_{4}
E_{\bar w}A_{5}E_{\bar\z}x\big\|^{2}\\
=&(\|y\|\prod_{j=1}^{3}\|A_{j}\|)^{2}\sum_{w,\z\in F}
\big\|E_{\z}A_{4}
E_{\bar w}A_{5}E_{\bar\z}x\big\|^{2}\\
\leq&(\|y\|\prod_{j=1}^{4}\|A_{j}\|)^{2}\sum_{\z\in F}\sum_{w\in F}
\big\|E_{\bar w}A_{5}E_{\bar\z}x\big\|^{2}\\
\leq&(\|y\|\prod_{j=1}^{5}\|A_{j}\|)^{2}\sum_{\z\in F}
\big\|E_{\bar\z}x\big\|^{2}
\leq(\|x\|\|y\|\prod_{j=1}^{5}\|A_{j}\|)^{2}\,,
\end{align*}
where the previous computation follows from the Schwarz and Bessel 
inequalities, and Pythagoras theorem. 
\end{proof}
\begin{lem}
\label{2}
The net $\big\{\sum_{z\in F}E_{z}AE_{\bar z}\,\big|\,F\,\text{finite 
subset of}\,\s_{\mathop{\rm pp}}^{\mathop{\rm a}}(U)\big\}$ 
converges in the strong operator 
topology.
\end{lem}
\begin{proof}
\begin{align*}
\big\|\sum_{z\in F}E_{z}AE_{\bar z}x-&\sum_{z\in G}E_{z}AE_{\bar z}x\big\|
\leq\big\|\sum_{z\in F\D G}E_{z}AE_{\bar z}x\big\|\\
+\big\|&\sum_{z\in F\backslash G}E_{z}AE_{\bar z}x\big\|
+\big\|\sum_{z\in G\backslash F}E_{z}AE_{\bar z}x\big\|\,.
\end{align*}

Thus, as the strong operator topology is complete, it is enough to prove that for
$\eps>0$, there exists a finite set $G_{\eps}$, such that
${\displaystyle\big\|\sum_{z\in H}E_{z}AE_{\bar z}x\big\|<\frac{\eps}{3}}$ whenever 
$H\subset G_{\eps}^{c}$. But,
$$
\big\|\sum_{z\in H}E_{z}AE_{\bar z}x\big\|^{2}
=\sum_{z\in H}\big\|E_{z}AE_{\bar z}x\big\|^{2}
\leq\|A\|^{2}\sum_{z\in H}\big\|E_{\bar z}x\big\|^{2}\,.
$$

The proof follows as the last sum is convergent. 
\end{proof}
\begin{prop}
\label{3}
The net 
$\big\{S^{F}_{\a;A_{1},\dots,A_{2k-1}}\,\big|\,F\,\text{finite 
subset of}\,\s_{\mathop{\rm pp}}^{\mathop{\rm a}}(U)\big\}$ converges in the weak operator 
topology for each finite set $\{A_{1},\dots,A_{2k-1}\}\subset\cb(\ch)$.
\end{prop}
\begin{proof}
By Lemma \ref{1} and Theorem II.1.3 of \cite{T}, it is enough to show that the 
$\big\{s^{F}_{A_{1},\dots,A_{2k-1}}(x,y)\big\}$ converges 
for each $x,y\in\ch$. We can also suppose without loss of generality that 
$x\in\ch$ is an eigenvector of $U$ with eigenvalue $z_{0}$.
The proof is by induction on $k$. By Lemma \ref{2}, it is enough to 
show that the assertion holds for the pair--partition 
$\b:\{1,\dots, 2k+2\}\mapsto\{1,\dots, k+1\}$, whenever it is true 
for any pair--partition $\a:\{1,\dots, 2k\}\mapsto\{1,\dots, k\}$. Let 
$k_{\b}\in\{1,\dots, 2k+2\}$ be the first element of the pair 
$\b^{-1}\big(\{k+1\}\big)$, and $\a_{\b}$ the pair--partition of
$\{1,\dots, 2k\}$ obtained 
by deleting $\b^{-1}\big(\{k+1\}\big)$ from $\{1,\dots, 2k+2\}$, and 
$k+1$ from $\{1,\dots, k+1\}$. We obtain
$$
s^{F}_{\b;A_{1},\dots,A_{2k+1}}(x,y)
=s^{F}_{\a_{\b};A_{1},\dots,A_{k_{\b}-1}
E_{\bar z_{0}}A_{k_{\b}+1},\dots,A_{2k}}(A_{2k+1}x,y)
$$
whenever $F$ is big enough, such that $\bar z_{0}\in F$. Thus in our 
situation, 
$$
\lim_{F\uparrow\s_{\mathop{\rm pp}}^{\mathop{\rm a}}(U)}
s^{F}_{\b;A_{1},\dots,A_{2k+1}}(x,y)=
\big\langle S_{\a_{\b};A_{1},\dots,A_{k_{\b}-1}
E_{\bar z_{0}}A_{k_{\b}-1},\dots,A_{2k}}A_{2k+1}x,y\big\rangle\,,
$$
$S_{\a;A_{1},\dots,A_{2k-1}}$ being the limit in the weak operator 
topology of $S^{F}_{\a;A_{1},\dots,A_{2k-1}}$ which exists by hypotesis.
\end{proof}

Proposition \ref{3} together with \eqref{opses}
allow us to define, and write symbolically for 
each finite subset 
$\{A_{1},\dots,A_{2k-1}\}\subset\ck(\ch)$,
\begin{align}
\label{symb}
&S_{\a;A_{1},\dots,A_{2k-1}}:=
\mathop{\rm w\!-\!lim}_{F\uparrow\s_{\mathop{\rm pp}}^{\mathop{\rm a}}(U)}
S^{F}_{\a;A_{1},\dots,A_{2k-1}}\\
=&\sum_{z_{1},\dots,z_{k}\in\s_{\mathop{\rm pp}}^{\mathop{\rm a}}(U)}
E_{z^{\#}_{\a(1)}}A_{1}E_{z^{\#}_{\a(2)}}\cdots
E_{z^{\#}_{\a(2k-1)}}A_{2k-1}E_{z^{\#}_{\a(2k)}}\,,\nn
\end{align}
where in \eqref{symb} the pairs $z^{\#}_{\a(i)}$ are alternatively $z_{j}$ and
$\bar z_{j}$ whenever $\a(i)=j$ as in \eqref{opses}.
\begin{prop}
\label{11}
If $A\in\ck(\ch)$, then
$$
\mathop{\rm w\!-\!lim}_{N}\frac{1}{N}\sum_{n=0}^{N-1}
U^{n}AU^{n}=S_{A}\,.
$$
\end{prop}
\begin{proof}
Consider ${\displaystyle\bigg\langle\frac{1}{N}\sum_{n=0}^{N-1}
U^{n}AU^{n}x,y\bigg\rangle}$.
By Lemma \ref{1}, we can suppose without loss of generality that 
$A=\langle\,\cdot\,,\eta\rangle\xi$. We can also suppose that some of 
the vectors $x,y,\xi,\eta$ are eigenvectors of $U$ if needed (see 
below). We have
$$
\bigg\langle\frac{1}{N}\sum_{n=0}^{N-1}
U^{n}AU^{n}x,y\bigg\rangle=\frac{1}{N}\sum_{n=0}^{N-1}
\langle U^{n}x,\eta\rangle\langle U^{n}\xi,y\rangle\,.
$$

We first suppose that $x,y,\xi,\eta\in\ch_{\mathop{\rm cont}}$, the last being the 
subspace of $\ch$ made of all vectors with continuous spectral 
measure on the unit circle (cf. \cite{RS}, Section VII.2). In this 
situation, we compute
\begin{align*}
&\frac{1}{N}\sum_{n=0}^{N-1}
\langle U^{n}x,\eta\rangle\langle U^{n}\xi,y\rangle\\
=\int\!\!\!\int_{\bt^{2}}&
\bigg(\frac{1}{N}\sum_{n=0}^{N-1}(zw)^{n}\bigg)
f_{x,\eta}(z)f_{\xi,y}(w)\di|\m_{x,\eta}|(z)
\di|\m_{\xi,y}|(w)\,,
\end{align*}
where $|\m_{x,\eta}|$, $|\m_{\xi,y}|$ are atomless positive 
bounded measures.\footnote{The measures $|\m_{x,y}|$, and the 
measurable functions $f_{x,y}$, $x,y\in\ch$ are 
the total variation measures of $\di\m_{x,y}(z):=\langle E(\di z)x,y\rangle$ 
and the corresponding densities, 
$\{E(z)\,|\,z\in\bt\}$ 
being the resolution of the identity of the unitary $U$.}  
As it was shown in \cite{L}, 
${\displaystyle\frac{1}{N}\sum_{n=0}^{N-1}(zw)^{n}}$ converges 
pointwise to indicator of the antidiagonal of the two--dimensional torus 
$\bt^{2}$, the last being negligible w.r.t. the product measure
$|\m_{x,\eta}|\times|\m_{\xi,y}|$. This means that if
$x,y,\xi,\eta\in\ch_{\mathop{\rm cont}}$, the ergodic mean under 
consideration is zero. The same happens if only one of the pairs $x,\eta$ 
or $\xi,y$ belongs to $\ch_{\mathop{\rm cont}}$. Namely, we suppose 
without loss of generality that $x,\eta\in\ch_{\mathop{\rm cont}}$ 
and, say $y$ is a eigenvector of $U$ with eigenvalue $w_{0}$. We have
\begin{align*}
&\frac{1}{N}\sum_{n=0}^{N-1}
\langle U^{n}x,\eta\rangle\langle U^{n}\xi,y\rangle\\
=\langle\xi,y\rangle\int_{\bt}&
\bigg(\frac{1}{N}\sum_{n=0}^{N-1}(zw_{0})^{n}\bigg)
f_{x,\eta}(z)\di|\m_{x,\eta}|(z)\,.
\end{align*}

In this situation, ${\displaystyle\frac{1}{N}\sum_{n=0}^{N-1}(zw_{0})^{n}}$ converges 
pointwise to the indicator of $\bar w_{0}$ which is negligible w.r.t. the measure
$|\m_{x,\eta}|$. Thus, we first conclude that 
$$
\frac{1}{N}\sum_{n=0}^{N-1}
\langle U^{n}x,\eta\rangle\langle U^{n}\xi,y\rangle=0
=\sum_{z\in\s_{\mathop{\rm pp}}^{\mathop{\rm a}}(U)}
\langle E_{z}x,\eta\rangle\langle E_{\bar z}\xi,y\rangle
\equiv\langle S_{A}x,y\rangle
$$
if at least one of the pairs $x,\eta$ 
or $\xi,y$ belongs to $\ch_{\mathop{\rm cont}}$. Second, 
${\displaystyle\frac{1}{N}\sum_{n=0}^{N-1}
\langle U^{n}x,\eta\rangle\langle U^{n}\xi,y\rangle}$
can be nonnull 
only if at least one of the elements of the pairs $x,\eta$ 
or $\xi,y$, say $\xi,\eta$, belong to $\ch_{\mathop{\rm pp}}$. As 
previously explained, we can suppose that $\xi,\eta$ are eigenvectors of $U$ with 
eigenvalues $z_{0}$, $w_{0}$. In this situation,   
\begin{align*}
&\frac{1}{N}\sum_{n=0}^{N-1}
\langle U^{n}x,\eta\rangle\langle U^{n}\xi,y\rangle=\langle x,\eta\rangle
\langle\xi,y\rangle\frac{1}{N}\sum_{n=0}^{N-1}(z_{0}w_{0})^{n}\\
\longrightarrow_{_{N}}&\langle x,\eta\rangle\langle\xi,y\rangle\d_{\bar z_{0},w_{0}}
=\sum_{z\in\s_{\mathop{\rm pp}}^{\mathop{\rm a}}(U)}
\langle E_{z}x,\eta\rangle\langle E_{\bar z}\xi,y\rangle
\equiv\langle S_{A}x,y\rangle\,.
\end{align*}

The proof follows by orthogonality, as the cases considered above exhaust all the 
possibilities.
\end{proof}

Here, there is the announced entangled ergodic theorem for 
pair--partitions and compact operators.
\begin{thm}
\label{main}
Let $\{A_{1},\dots,A_{2k-1}\}\subset\ck(\ch)$, and 
$\a:\{1,\dots, 2k\}\mapsto\{1,\dots, k\}$ a pair--partition of the 
set $\{1,\dots, 2k\}$. Then
\begin{align*}
\mathop{\rm w\!-\!lim}_{N}&\bigg\{\frac{1}{N^{k}}\sum_{n_{1},\dots,n_{k}=0}^{N-1}
U^{n_{\a(1)}}A_{1}U^{n_{\a(2)}}\cdots 
U^{n_{\a(2k-1)}}A_{2k-1}U^{n_{\a(2k)}}\bigg\}\\
=&S_{\a;A_{1},\dots,A_{2k-1}}\,,
\end{align*} 
where $S_{\a;A_{1},\dots,A_{2k-1}}$ 
is given in \eqref{symb}.
\end{thm}
\begin{proof}
We start by noticing that 
$$
\bigg\|\frac{1}{N^{k}}\sum_{n_{1},\dots,n_{k}=0}^{N-1}
U^{n_{\a(1)}}A_{1}U^{n_{\a(2)}}\cdots 
U^{n_{\a(2k-1)}}A_{2k-1}U^{n_{\a(2k)}}\bigg\|\leq\prod_{j=1}^{2k-1}\|A_{j}\|\,.
$$

Thanks to this and Lemma \ref{1}, 
as $\{A_{1},\dots,A_{2k-1}\}\subset\ck(\ch)$, we can suppose 
that the $A_{j}$ are rank one. 
Indeed, put $K:=\max_{1\leq j\leq 2k-1}\|A_{j}\|$. Choose finite rank 
operators $A^{\eps}_{j}$, such that $\|A^{\eps}_{j}\|\leq K$ and
$\|A_{j}-A^{\eps}_{j}\|<{\displaystyle\frac{\eps}{4(2k-1)K^{2k}\|x\|\|y\|}}$,
$j=1,\dots,2k-1$. Then
\begin{align*}
&\bigg|\bigg\langle\frac{1}{N^{k}}\sum_{n_{1},\dots,n_{k}=0}^{N-1}
U^{n_{\a(1)}}A_{1}U^{n_{\a(2)}}\cdots 
U^{n_{\a(2k-1)}}A_{2k-1}U^{n_{\a(2k)}}x,y\bigg\rangle\\
-&\big\langle S_{\a;A_{1},\dots,A_{2k-1}}x,y\big\rangle\bigg|\leq\frac{\eps}{2}\\
+&\bigg|\bigg\langle\frac{1}{N^{k}}\sum_{n_{1},\dots,n_{k}=0}^{N-1}
U^{n_{\a(1)}}A^{\eps}_{1}U^{n_{\a(2)}}\cdots 
U^{n_{\a(2k-1)}}A^{\eps}_{2k-1}U^{n_{\a(2k)}}x,y\bigg\rangle\\
-&\big\langle S_{\a;A^{\eps}_{1},\dots,A^{\eps}_{2k-1}}x,y\big\rangle\bigg|\,.
\end{align*}

Thus, for rank one operators $A_{j}$, we obtain
\begin{align*}
&\bigg\langle\frac{1}{N^{k}}\sum_{n_{1},\dots,n_{k}=0}^{N-1}
U^{n_{\a(1)}}A_{1}U^{n_{\a(2)}}\cdots 
U^{n_{\a(2k-1)}}A_{2k-1}U^{n_{\a(2k)}}x,y\bigg\rangle\\
=&\prod_{j=1}^{k}
\frac{1}{N}\sum_{n_{j}=0}^{N-1}\big\langle U^{n_{j}}x_{j},y_{j}\big\rangle
\big\langle U^{n_{j}}\xi_{j},\eta_{j}\big\rangle
\equiv\prod_{j=1}^{k}\frac{1}{N}\sum_{n_{j}=0}^{N-1}
\big\langle U^{n_{j}}B_{j}U^{n_{j}}\xi_{j},y_{j}\big\rangle
\end{align*}
with $B_{j}=\langle\,\cdot\,,\eta_{j}\rangle x_{j}$. 
Here, $\{x_{j},y_{j},\xi_{j},\eta_{j}\}_{j=1}^{k}$
are suitable vectors depending on the $A_{i}$ and $x$, $y$.
By Proposition \ref{11}, we obtain
\begin{align*}
&\bigg\langle\frac{1}{N^{k}}\sum_{n_{1},\dots,n_{k}=0}^{N-1}
U^{n_{\a(1)}}A_{1}U^{n_{\a(2)}}\cdots 
U^{n_{\a(2k-1)}}A_{2k-1}U^{n_{\a(2k)}}x,y\bigg\rangle\\
=&\prod_{j=1}^{k}\frac{1}{N}\sum_{n_{j}=0}^{N-1}
\big\langle U^{n_{j}}B_{j}U^{n_{j}}\xi_{j},y_{j}\big\rangle
\longrightarrow_{_{N}}\prod_{j=1}^{k}\big\langle S_{B_{j}}\xi_{j},y_{j}\big\rangle\\
=&\big\langle S_{\a;A_{1},\dots,A_{2k-1}}x,y\big\rangle\,.
\end{align*}
\end{proof}

We end the present section by proving the entangled ergodic theorem 
for general partitions of any finite set $\{1,\dots,m\}$, and for 
compact operators $\{A_{1},\dots,A_{m-1}\}$.
\begin{thm}
\label{main1}
For $m\geq k$, let $\a:\{1,\dots,m\}\mapsto\{1,\dots,k\}$ be a partition of the 
set $\{1,\dots,m\}$. If $\{A_{1},\dots,A_{m-1}\}\subset\ck(\ch)$, then the
ergodic average 
\begin{equation}
\label{wl}
\frac{1}{N^{k}}\sum_{n_{1},\dots,n_{k}=0}^{N-1}
U^{n_{\a(1)}}A_{1}U^{n_{\a(2)}}\cdots 
U^{n_{\a(m-1)}}A_{m-1}U^{n_{\a(m)}}
\end{equation}
converges in the weak operator topology.
\end{thm}
\begin{proof}
As before, it is enough to show that  
$$
\bigg\langle\frac{1}{N^{k}}\sum_{n_{1},\dots,n_{k}=0}^{N-1}
U^{n_{\a(1)}}A_{1}U^{n_{\a(2)}}\cdots 
U^{n_{\a(m-1)}}A_{m-1}U^{n_{\a(m)}}x,y\bigg\rangle
$$
converges for every $x,y\in\ch$, whenever the $A_{j}$ are rank one
operators. But, in this situation,
\begin{align*}
&\bigg\langle\frac{1}{N^{k}}\sum_{n_{1},\dots,n_{k}=0}^{N-1}
U^{n_{\a(1)}}A_{1}U^{n_{\a(2)}}\cdots 
U^{n_{\a(m-1)}}A_{m-1}U^{n_{\a(m)}}x,y\bigg\rangle\\
=&\prod_{j=1}^{k}\frac{1}{N}\sum_{n_{j}=0}^{N-1}
\prod_{\{p\,|\,\a(p)=j\}}\big\langle U^{n_{j}}x_{p,j},y_{p,j}\big\rangle\\
=&\prod_{j=1}^{k}\iint\cdot\cdot\int_{\bt^{|\a^{-1}\{j\}|}}
\bigg(\frac{1}{N}\sum_{n_{j}=0}^{N-1}
\big(\prod_{\{p\,|\,\a(p)=j\}}z_{p}\big)^{n_{j}}\bigg)
\prod_{\{p\,|\,\a(p)=j\}}\langle E(\di z_{p})x_{p,j},y_{p,j}\rangle\\
&\longrightarrow_{{}_{N}}
\prod_{j=1}^{k}\iint\cdot\cdot\int_{\bt^{|\a^{-1}\{j\}|}}
\chi_{\{1\}}\bigg(\prod_{\{p\,|\,\a(p)=j\}}z_{p}\bigg)
\prod_{\{p\,|\,\a(p)=j\}}\langle E(\di z_{p})x_{p,j},y_{p,j}\rangle
\end{align*}
where we have 
used the Lebesgue dominated convergence theorem. 
Here, the $x_{p,j}$, $y_{p,j}$ are vectors uniquely determined by the 
rank one operators $A_{1},\dots,A_{m-1}$ and vectors $x,y$, and 
$\chi_{\G}$ denotes the indicator of the set $\G$. 
\end{proof}

We notice that it seems difficult to provide
an esplicit formula for the weak limit of \eqref{wl} similar to that 
in Theorem \ref{main} relative to the case of 
pair--partitions.

\section{the almost periodic case}

The present section deals with some cases relative to the almost 
periodic situation, without any restriction relative to the operators 
appearing in \eqref{0}. We then suppose that $\ch$ is generated by 
the eigenvectors of $U$.
\begin{prop}
\label{qper2}    
In the almost periodic case,
$$
\mathop{\rm s\!-\!lim}_{N}\frac{1}{N}\sum_{n=0}^{N-1}
U^{n}AU^{n}=S_{A}
$$
for each $A\in\cb(\ch)$.
\end{prop}
\begin{proof}
By our assumptions, we can suppose that $x$ is an eigenvector of $U$ 
with eigenvalue $z_{0}$. We have
\begin{align*}
&\frac{1}{N}\sum_{n=0}^{N-1}U^{n}AU^{n}x
=\frac{1}{N}\sum_{n=0}^{N-1}(z_{0}U)^{n}Ax\\
=&\bigg(\frac{1}{N}\sum_{n=0}^{N-1}(z_{0}U)^{n}\bigg)Ax
\longrightarrow_{{}_{N}} E_{\bar z_{0}}Ax\equiv S_{A}x\,.
\end{align*}

Here, we have used the mean ergodic theorem (cf. \cite{RS}), and as 
usual, $E_{z}\equiv E(\{z\})$ is the selfadjoint projection onto the 
eigenspace corresponding to the eigenvalue $z\in\bt$.
\end{proof}

Now we treat the cases relative 
to all the pair partitions of four elements.
\begin{thm}
\label{qper}
Suppose that the 
dynamics induced by the unitary $U$ on $\ch$ is almost periodic. Then 
\begin{align*}
\mathop{\rm w\!-\!lim}_{N}&\bigg\{\frac{1}{N^{2}}\sum_{n_{1},n_{2}=0}^{N-1}
U^{n_{\a(1)}}AU^{n_{\a(2)}}B 
U^{n_{\a(3)}}CU^{n_{\a(4)}}\bigg\}\\
=&S_{\a;A,B,C}
\end{align*} 
for each pair--partition $\a$ of four elements, and every 
$\{A,B,C\}\subset\cb(\ch)$.
\end{thm}
\begin{proof}
As previously explained, we can suppose that $x,y\in\ch$ are
eigenvectors of $U$ with eigenvalues $z_{0},w_{0}$, respectively.

Suppose that $\a$ is the partition $\{1,1,2,2\}$. Then, by 
Proposition \ref{qper2},
\begin{align*}
&\bigg\langle\frac{1}{N^{2}}\sum_{k,n=0}^{N-1}
U^{k}AU^{k}B 
U^{n}CU^{n}x,y\bigg\rangle\\
=&\bigg\langle\bigg(\frac{1}{N}\sum_{k=0}^{N-1}
U^{k}AU^{k}\bigg)B 
\bigg(\frac{1}{N}\sum_{n=0}^{N-1}
U^{n}CU^{n}\bigg)x,y\bigg\rangle\\
&\longrightarrow_{{}_{N}}\big\langle S_{A}BS_{C}x,y\big\rangle
\equiv\big\langle S_{\a;A,B,C}x,y\big\rangle\,.
\end{align*}

Let $\a$ be the partition $\{1,2,2,1\}$. Then, again by Proposition 
\ref{qper2},
\begin{align*}
&\bigg\langle\frac{1}{N^{2}}\sum_{k,n=0}^{N-1}
U^{k}AU^{n}B 
U^{n}CU^{k}x,y\bigg\rangle\\
=&\bigg(\frac{1}{N}\sum_{k=0}^{N-1}
(z_{0}w_{0})^{k}\bigg)
\bigg\langle\frac{1}{N}\sum_{n=0}^{N-1}
AU^{n}BU^{n}Cx,y\bigg\rangle\\
&\longrightarrow_{{}_{N}}\d_{\bar z_{0},w_{0}}\big\langle AS_{B}Cx,y\big\rangle
\equiv\big\langle S_{\a;A,B,C}x,y\big\rangle\,.
\end{align*}

Finally, if $\a$ is the partition $\{1,2,1,2\}$, then by the mean ergodic 
theorem,
\begin{align*}
&\bigg\langle\frac{1}{N^{2}}\sum_{k,n=0}^{N-1}
U^{k}AU^{n}B 
U^{k}CU^{n}x,y\bigg\rangle\\
=&\bigg\langle A\bigg(\frac{1}{N}\sum_{k=0}^{N-1}(z_{0}U)^{k}\bigg)B
\bigg(\frac{1}{N}\sum_{n=0}^{N-1}(w_{0}U)^{n}\bigg)Cx,y\big\rangle\\
&\longrightarrow_{{}_{N}}\big\langle AE_{\bar z_{0}}BE_{\bar w_{0}}Cx,y\big\rangle
\equiv\big\langle S_{\a;A,B,C}x,y\big\rangle\,.
\end{align*}
\end{proof}

To end the present section by noticing that a general entangled 
ergodic theorem is not yet available even in the almost periodic case. 
However, Proposition \ref{qper2} and Theorem \ref{qper} allow us to 
treat, always in the almost periodic case, other situations relative 
to pair--partitions of sets with more than four elements.

\section{outlook}

\end{document}